\documentclass[a4paper,11pt]{amsart}

\usepackage{amsmath,amssymb,amsfonts}
\usepackage{amscd}
\usepackage[all,cmtip]{xy}
\usepackage{framed}
\usepackage[inline]{enumitem}  
\usepackage{float}
\usepackage{multicol}
\usepackage{array}
\usepackage{makecell}
\newcolumntype{P}[1]{>{\centering\arraybackslash}p{#1}}
\newcolumntype{M}[1]{>{\centering\arraybackslash}m{#1}}
\usepackage{lscape}
\usepackage{times}
\usepackage{xcolor}
\usepackage[verbose]{hyperref}
\hypersetup{colorlinks=false,allbordercolors=blue,pdfborderstyle={/S/U/W 1}}
\usepackage{setspace}

\theoremstyle{plain}
\newtheorem{thm}{Theorem}[section]
\newtheorem{theorem}[thm]{Theorem}

\newtheorem{lemma}[thm]{Lemma}
\newtheorem{corollary}[thm]{Corollary}
\newtheorem{proposition}[thm]{Proposition}

\theoremstyle{definition}
\newtheorem*{notation*}{Notation}
\newtheorem*{notation}{Notation}
\newtheorem{remark}[thm]{Remark}

\newtheorem*{claim*}{Claim}
\newtheorem*{claimproof*}{Proof of Claim}

\newtheorem{example}[thm]{Example}

\numberwithin{equation}{section}
\newcommand{\C}{{\mathbb C}}
\newcommand{\Z}{{\mathbb Z}}
\newcommand{\rH}{{\rm Hom}}
\newcommand{\id}{\operatorname{id}}
\newcommand{\tr}{\operatorname{trace}}
\newcommand{\transpose}{\operatorname{Tr}}

\begin{document}

\title{Matrix representations of Lucas sequences \\and Salem traces on K3 surfaces}
\author{Kwangwoo Lee}
\address{Center for Complex Geometry,
Institute for Basic Science,
55 Expo-ro, Yuseong-gu, Daejeon 34126,
Republic of Korea}
\email{klee0222@gmail.com}


%\date{\today}

\maketitle

\begin{abstract}
We develop matrix representations of Lucas sequences with arbitrary initial terms and investigate their connections to associated Pell-type equations. Using these representations, we reformulate Pell-type equations in terms of trace and determinant expressions, yielding a generalization of a classical characterization of Fibonacci numbers. We then study solutions to such Pell-type equations via Lucas sequences with arbitrary initial terms. We also construct K3 surfaces admitting automorphisms induced by these matrix representations and describe the associated Salem traces explicitly in terms of Lucas sequences.
\end{abstract}

\section{Introduction}\label{introduction}
We consider the Diophantine equation
\begin{equation}\label{Pell eq}
v^2-Du^2=\pm 4N,
\end{equation}
where $D$ and $N$ are positive integers, and $D$ is not a perfect square. This equation naturally arises in the study of Pell-type equations, binary quadratic forms, and arithmetic problems on algebraic varieties.

Let $(u,v)$ be a solution of \eqref{Pell eq}. Multiplying the corresponding algebraic integer $1/2(v+u\sqrt{D})$ by any solution $1/2(y+x\sqrt{D})$ of the Pell equation $y^2-Dx^2=4$ produces another solution of \eqref{Pell eq}. 
Hence, the set of solutions is partitioned into equivalence classes under
the action of the unit group. Each class is generated from its unique
fundamental solution, namely the solution whose $u$-value is
the smallest non-negative among all solutions in the class. In \cite{St}, Stolt found some conditions on the fundamental solutions, hence such fundamental solutions can be found after a finite number of trials.

In this paper, we characterize the fundamental solutions in terms of Lucas pairs. Each Lucas pair may represent a fundamental solution of a solution class, and all other solutions in the class are generated recursively from it. This avoids the action of the unit group and provides a simpler and more efficient description of the complete solution set. Moreover, we construct K3 surfaces admitting automorphisms whose associated Salem traces are given in terms of Lucas sequences.

For $p,q \in \Z\setminus{0}$ and $a,b\in \Z$, we have two Lucas sequences $
\{u_n^{a,b}(p,q)\}_{n\in \Z}$ and $\{v_n^{a,b}(p,q)\}_{n\in \Z}$ defined by
\begin{equation}\label{Lucas seq}
\begin{gathered}
u_0^{a,b} = a, ~\ u_1^{a,b} = b, ~\ u_{n+1}^{a,b} = p u_n^{a,b} - q u_{n-1}^{a,b}, ~\ n\in \Z; \\
v_0^{a,b} = 2b - ap, ~\ v_1^{a,b} = bp - 2aq, ~\ v_{n+1}^{a,b} = p v_n^{a,b} - q v_{n-1}^{a,b}, ~\ n\in \Z.
\end{gathered}
\end{equation}
If $(u_n^{a,b}(p,q),v_n^{a,b}(p,q))$ is a solution of \eqref{Pell eq} for some $n\in \Z$, then for all integers $n$, such pairs are solutions (cf. Lemma~\ref{main theorem 2}). Hence, a pair $(a,b)$ of initial terms defines a class of solutions $\{(u_n^{a,b}(p,q),v_n^{a,b}(p,q))\}_{n\in \Z}$.
Moreover, we define an equivalence relation of two pairs $(a,b)$ and $(a',b')$ denoted by $(a,b)\sim (a',b')$ according to Lemma~\ref{sim}.

To state our main theorems, fix $p \in \Z_{>0}$ and $q \in \{\pm 1\}$, and set $D = p^2 - 4q$. For a given positive integer $N$, define
\begin{gather}
T:=\{(a,b)\in \Z^2 \mid  ~ (2b-ap)^2-a^2D=4N \}, \label{T} \\
S:=\{(x,y)\in \Z^2 \mid ~ y^2-Dx^2=4 N q^n \text{ for some } n\in \Z \}.
\end{gather}
We first establish the following lemma.
\begin{lemma}\label{main theorem 2} For positive integers $p$ and $N$, and $q \in \{\pm 1\}$,
\begin{equation}
S=\bigcup_{(a,b)\in T} \{ (\epsilon_1 u_n^{a,b}(p,q), \epsilon_2 v_n^{a,b}(p,q)) \mid n \in \mathbb{Z},\ \epsilon_1, \epsilon_2 \in \{\pm 1\} \}
\end{equation}
is the set of all solutions of \eqref{Pell eq} with $N=b^2-abp+a^2q$.
\end{lemma}

Thus, all solutions of \eqref{Pell eq} arise from Lucas sequences with suitable initial terms. Unlike the case $(a,b)=(0,1)$, however, there are in general many fundamental solutions.
% Although Theorem~\ref{main theorem 2} describes all solutions, it does not immediately yield a finite and effective procedure.
We therefore develop explicit methods to compute a finite set of representatives as follows.

\begin{theorem}\label{method}
Assume that $q=-1$, and let $(a_0,b_0)\in T$
be such that $a_0\ge0$ is minimal and 
$0\le a_0<b_0$. %, and $b_0^2-pa_0b_0-a_0^2=N$.
Set
\begin{equation}\label{condition T'}
T':=
\{(a,b)\in T\mid
a_0\le a<b_0
\text{ or }
a_0\le b<b_0\},
\end{equation}
and define $\sim$ on $T'$ as in Lemma~\ref{sim}. 
Then
\begin{equation}\label{simple form of solutions 1}
S=
\bigsqcup_{(a,b)\in T'/\sim}
\{
(\epsilon_1u_n^{a,b}(p,-1),
\epsilon_2v_n^{a,b}(p,-1))
\mid
n\in\Z,\;
\epsilon_1,\epsilon_2\in\{\pm1\}
\}.
\end{equation}
\end{theorem}

\begin{theorem}\label{method2}
Assume that $q=1$, and let $(a_0,b_0)\in T$ be such that $a_0\ge0$ is minimal.
Set
\begin{equation}\label{alpha}
T'':=\{(a,b)\in T \mid a_0\le a<b_0\},
\end{equation}
and define $\sim$ on $T''$ as in Lemma~\ref{sim}. 
Then 
\begin{equation}\label{simple form of solutions 2}
S=
\bigsqcup_{(a,b)\in T''/\sim}
\left\{
(\epsilon_1u_n^{a,b}(p,1),\epsilon_2v_n^{a,b}(p,1))
\mid
n\in\mathbb Z,\;
\epsilon_1,\epsilon_2\in\{\pm1\}
\right\}.
\end{equation}
\end{theorem}

These are generalizations of Theorem~\ref{solutions to Pell's equation in Lucas seq} or the representation of solutions \eqref{S_0}.
Moreover, using the matrix representations of Lucas sequences, we have the following theorem.
\begin{theorem}(Theorem \ref{K3 automorphisms from matrix representations})\label{main thm 3}
Suppose that $q \neq \pm 1$, $q \mid p$, $D=p^2-4q>0$, and assume that $D$ is not a square number. Then there exists a K3 surface $X$ with a symplectic automorphism $g$ whose Salem trace is $(\frac{v_{2n}^{0,1}}{q^n})^2-2$.
\end{theorem}

The paper is organized as follows. Sections~\ref{Lucas Sequences and Their Matrix Representations} and \ref{Pell's equations} review Lucas sequences, their matrix representations, and Pell equations. Moreover, we generalize a classical characterization of Fibonacci numbers to Lucas numbers with arbitrary initial terms (cf. Corollary~\ref{Whitney theorem}). Section~\ref{Proofs} proves Lemma~\ref{main theorem 2}, Theorems~\ref{method},~\ref{method2}, and provides examples. Finally, in Section~\ref{Automorphisms of K3 Surfaces}, we prove Theorem~\ref{main thm 3}.

\begin{notation}\label{notation}
Throughout this paper, we assume that $D=p^2-4q$ is positive and not a square. We use the following notation:
\begin{itemize}
\item $u_n^{a,b}:=u_n^{a,b}(p,q)$, ~ $v_n^{a,b}:=v_n^{a,b}(p,q)$,  ~ $M:=M_{(p,q)}$;
\item $u_n:=u_n^{0,1}(p,q)$, ~ $v_n:=v_n^{0,1}(p,q)$, ~ $f_n:=u_n^{0,1}(1,-1)$, ~ $l_n:=v_n^{0,1}(1,-1)$.
\end{itemize}
\end{notation}

\section{Matrix representations of Lucas sequences}\label{Lucas Sequences and Their Matrix Representations}

In \cite{S}, Silvester first showed many properties of the Fibonacci numbers $f_n$ using a matrix representation. Indeed, for  $M_{(1,\,-1)}=\begin{pmatrix} 0& ~1\\ 1&~1\end{pmatrix}$,
he showed that
\begin{equation}
M_{(1,\,-1)}^n= \begin{pmatrix}
   f_{n-1}      & ~f_n    \\
   f_n    & ~f_{n+1}    
  \end{pmatrix}, \: n\geq 1.
\end{equation}
For any integers $p$ and $q$, the sequence $\{u_n\}_{n \ge 0}$ with $u_0 = 0$ and $u_1 = 1$ has the following matrix representation, analogous to that of the Fibonacci sequence.

Let $M := M_{(p,q)} = \begin{pmatrix} 0 & 1 \\ -q & p \end{pmatrix}$. Then, for $n \geq 1$,
\begin{equation}\label{power of M}
M^n =
\begin{pmatrix}
 -q u_{n-1} & u_n \\
 -q u_n & u_{n+1}
\end{pmatrix}.
\end{equation}
In fact, this identity remains valid for all $n \in \mathbb{Z}$.

We extend \eqref{power of M} to Lucas sequences with arbitrary initial terms. By the recurrence \eqref{Lucas seq}, $-qu_{-1}^{a,b}=b-ap$ for initial terms $u_0^{a,b}=a$ and $u_1^{a,b}=b$. We define
\begin{equation}\label{E}
E_{(a,b)}:=\begin{pmatrix}
   -qu_{-1}^{a,b}  &~ u_0^{a,b} \\
   -qu_0^{a,b} &~ u_1^{a,b} 
  \end{pmatrix}=
   \begin{pmatrix}
   b-ap  &~ a \\
   -qa &~ b 
  \end{pmatrix}.
\end{equation}

\begin{lemma}\label{power of M and E}
For $\{u_n^{a,b}\}_{n\in \Z}$, we have 
\begin{equation}\label{equation of power of M and E}
M^nE_{(a,b)}=\begin{pmatrix}
   -q u_{n-1}^{a,b}  &~ u_n^{a,b} \\
   -q u_n^{a,b} &~ u_{n+1}^{a,b}
  \end{pmatrix},  ~\ n\in \Z.
\end{equation}
\end{lemma}

\begin{proof}
The result follows from the recurrence relation for $\{u_n^{a,b}\}$ in \eqref{Lucas seq}.
\end{proof}

Similarly, the sequence $\{v_n^{a,b}\}_{n \in \Z}$ also admits a matrix representation. 
Note that $\{v_n^{a,b}\}_{n \in \Z}$ satisfies the same recurrence relation as $\{u_n^{a,b}\}_{n \in \Z}$, but with different initial terms $v_0^{a,b}=2b-ap, v_1^{a,b}=bp-2aq$.
Hence we obtain the following result.

\begin{lemma}\label{matrix representation for v_n}
For $\{v_n^{a,b}\}_{n\in \Z}$, we have 
\begin{equation}\label{equation of matrix representation for v_n}
M^nE'_{(a,b)}=\begin{pmatrix}
   -qv_{n-1}^{a,b}  &~ v_n^{a,b} \\
   -qv_n^{a,b} &~ v_{n+1}^{a,b}
  \end{pmatrix}, ~\ n\in \Z, 
\end{equation}
where 
\begin{equation}\label{E'}
E'_{(a,b)}:=\begin{pmatrix}
   -qv_{-1}^{a,b}  &~ v_0^{a,b} \\
   -qv_0^{a,b} &~ v_1^{a,b} 
  \end{pmatrix}=\begin{pmatrix}
   -bp+ap^2-2aq  &~ 2b-ap \\
   -2bq+apq &~ bp-2aq 
  \end{pmatrix}.
\end{equation}
\end{lemma}

Using the matrix representations of Lucas sequences, we obtain the following properties.

\begin{lemma}\label{properties of Lucas seq1}
For Lucas sequences $\{u_n^{a,b}\}_{n \in \Z}$ and $\{v_n^{a,b}\}_{n \in \Z}$ associated with $p, q \in \Z\setminus 0$, the following properties hold:
\begin{enumerate}
\item $v^{a,b}_n=u^{a,b}_{n+1}-qu^{a,b}_{n-1}=\tr(M^nE_{(a,b)})$;
\item $\tr((M^nE_{(a,b)})^2)=D(u_n^{a,b})^2+2\det(M^nE_{(a,b)})$.
\end{enumerate}
\end{lemma}

\begin{proof}
(1). From the matrix representation in Lemma~\ref{matrix representation for v_n}, we have that $v^{a,b}_n$ is the $(1,2)$-entry of $M^n E'_{(0,1)}$. By multiplying $M^n$ and $E'_{(0,1)}$, we obtain
\begin{equation*}
v^{a,b}_n=-qu^{0,1}_{n-1}v^{a,b}_0+u^{0,1}_nv^{a,b}_1=-qu^{0,1}_{n-1}(2b-ap)+u^{0,1}_n(bp-2aq)
\end{equation*}
which is the trace of $M^nE_{(a,b)}$. 

\noindent(2). From the identity \eqref{equation of power of M and E}, we have
\begin{align}\label{multiequation in useful lemma}
\tr((M^nE_{(a,b)})^2)=&q^2(u_{n-1}^{a,b})^2-2q(u_{n}^{a,b})^2+(u_{n+1}^{a,b})^2\nonumber\\
=&q^2(u_{n-1}^{a,b})^2-2q(u_{n}^{a,b})^2+(pu_n^{a,b}-qu_{n-1}^{a,b})^2\nonumber\\
=&(p^2-2q)(u_n^{a,b})^2-2qu_{n-1}^{a,b}u_{n+1}^{a,b}.
\end{align}
On the other hand, by taking the determinant on both sides of \eqref{equation of power of M and E}, we obtain 
\[qu_{n-1}^{a,b}u_{n+1}^{a,b}=-\det(M^nE_{(a,b)})+q(u_{n}^{a,b})^2.\] 
Substituting this back into \eqref{multiequation in useful lemma}, we obtain the desired result.
\end{proof}

\begin{lemma}\label{sim}
Let $n\in \Z$ and $q, \epsilon_1, \epsilon_2\in \{\pm 1\}$. For two pairs $(a,b)$ and $(a',b')$ in $T$, the following are equivalent:
\begin{enumerate}
\item $\{(\epsilon_1 u_n^{a,b}, \epsilon_2 v_n^{a,b})\}=\{(\epsilon_1 u_n^{a',b'}, \epsilon_2 v_n^{a',b'})\}$;
\item $\{(\epsilon_1 u_n^{a,b}, \epsilon_2 v_n^{a,b})\}\cap \{(\epsilon_1 u_n^{a',b'}, \epsilon_2 v_n^{a',b'})\}\neq \emptyset$;
\item $(a',b')=\pm(u_k^{a,b}, u_{k+1}^{a,b})$ ~\ or ~\ $\pm(u_k^{a,b}, qu_{k-1}^{a,b})$ ~\ for some $k\in \Z$.
\end{enumerate}
If one of the above conditions holds, we denote $(a,b)\sim (a',b')$.
\end{lemma}

\begin{proof}
(1)$\Rightarrow$(2) and (3)$\Rightarrow$(1) are trivial. 

For (2)$\Rightarrow$(3), suppose that $(u_m^{a,b},v_m^{a,b})=(u_n^{a',b'},v_n^{a',b'})$ for some $m,n\in \Z$ up to sign. %Since they are solutions of \eqref{Pell eq}, we have that $Q(a,b)=\pm Q(a',b')$ by Lemma~\ref{main theorem 2}. 
By Lemma~\ref{properties of Lucas seq1}(1), $v_m^{a,b}=u_{m+1}^{a,b}-qu_{m-1}^{a,b}$, and by \eqref{Lucas seq}, \begin{equation}
v_m^{a,b}=pu_{m}^{a,b}-2qu_{m-1}^{a,b}~\ \text{ and } ~\ v_n^{a',b'}=pu_{n}^{a',b'}-2qu_{n-1}^{a',b'}.
\end{equation} 

Now, if $(u_m^{a,b},v_m^{a,b})=\pm (u_n^{a',b'},v_n^{a',b'})$, we have $u_{m-1}^{a,b}=\pm u_{n-1}^{a',b'}$. Hence, by shifting the index, we have $(a',b')=\pm (u_k^{a,b}, u_{k+1}^{a,b})$ for some $k$. If $(u_m^{a,b},v_m^{a,b})=\pm(u_n^{a',b'},-v_n^{a',b'})$, we have $u_{m+1}^{a,b}=\pm qu_{n-1}^{a',b'}$. Hence, by shifting the index, we have $(a',b')=\pm (u_k^{a,b}, qu_{k-1}^{a,b})$ for some $k$. 
\end{proof}

In \cite{W}, Whitney showed that a given integer $n$ is a Fibonacci number if and only if either $5n^2+4$ or $5n^2-4$ is a square number. We extend this characterization to Lucas numbers as follows.

\begin{corollary}\label{Whitney theorem}
An integer $m$ satisfies
\begin{equation}\label{equation in Whitney theorem}
Dm^2+4\det(M^nE_{(a,b)})=(\tr M^nE_{(a,b)})^2
\end{equation}
if and only if $m=u_n^{a,b}$.
\end{corollary}

\begin{proof}
This follows from Lemma~\ref{properties of Lucas seq1}(2) together with the identity
\begin{equation}\label{matrix identity}
(\tr A)^2-\tr(A^2)=2\det (A),
\end{equation} 
applied to $A=M^nE_{(a,b)}$. Moreover, note that $v_n^{a,b}=\tr(M^nE_{(a,b)})$ by Lemma \ref{properties of Lucas seq1}(1).
\end{proof}

\section{Diophantine equations of Pell type}\label{Pell's equations}

Let $D$ be any positive integer. The Pell equations associated to $D$ are
\begin{align}
v^2-Du^2=&4 \hspace{14mm}  \text{ the positive Pell equation};\label{positive Pell equation}\\
v^2-Du^2=&-4   \hspace{10mm}\text{ the negative Pell equation}.\label{negative Pell equation}
\end{align}
The solutions are described in \cite[Theorem 3.18]{B}. 
In particular, when $D=p^2-4q$ for some $p>0$ and $q=\pm 1$, the solutions can be written in terms of Lucas sequences $\{u_n\}$ and $\{v_n\}$ as follows.

\begin{theorem} [{\cite[Theorem 1 and Theorem 7]{A}}]\label{solutions to Pell's equation in Lucas seq}
Let $p, q$ be integers such that $p>0, \vert q\vert =1, (p,q)\neq (3,1)$, and $D=p^2-4q>0$. If positive integers $u$ and $v$ are such that 
\begin{equation}
v^2-Du^2=\pm 4,
\end{equation}
then 
\begin{equation}
u\in \{u_n\}_{n\geq 0}~\ \text{ and }~\ v\in \{v_n\}_{n\geq 0}.
\end{equation}
Moreover, $(u_n(3,1),v_n(3,1))=(u_{2n}(1,-1),v_{2n}(1,-1))$.
\end{theorem}

\begin{remark}\label{remark for solutions 4} For $p>0$, $|q|=1$, $(p,q)\neq (3,1)$, and $(a,b)=(0,1)$,
Theorem~\ref{solutions to Pell's equation in Lucas seq} implies that
\begin{equation}\label{S_0}
\{(\epsilon_1u_n, \epsilon_2v_n)\mid n\in \Z, \epsilon_1, \epsilon_2\in \{\pm 1\} \}=S_0
\end{equation} 
is the set of all solutions of $y^2-(p^2-4q)x^2=4q^n$.
\end{remark}

\section{Solutions via Lucas sequences}\label{Proofs}

We obtain a similar description for the generalized Pell-type equation
\begin{equation}\label{gen Pell's equ}
y^2-Dx^2=4N q^n,
\end{equation}
where $N=b^2-abp+a^2q$ for some $a,b\in \Z$. 
We need the following lemma for the proof of Lemma~\ref{main theorem 2}.

\begin{lemma}\label{description of solution set}
Let $(x_0,y_0)$ be an integer solution to \eqref{gen Pell's equ}
with $q=\pm1$.
Then there exists a pair $(a',b') \in T$, as defined in the introduction, such that 
\begin{equation}
(x_0,y_0) = (u_0^{a',b'}, v_0^{a',b'}) \quad \text{or} \quad (u_1^{a',b'}, v_1^{a',b'}).
\end{equation}
\end{lemma}

\begin{proof}
We consider the cases $q = \pm 1$ separately. Suppose that $(x_0, y_0)$ is a solution to \eqref{gen Pell's equ}.

%\medskip
\noindent\textbf{Case 1:} $q = 1$. 
Then $(x_0,y_0)$ satisfies
\[
y_0^2-(p^2-4)x_0^2=4N.
\]
We define $a':=x_0$ and determine $b'$ by
$y_0=2b'-a'p$, from \eqref{Lucas seq},
which gives
\[
b'=\frac{p x_0+y_0}{2}.
\]
Since $y_0^2 \equiv p^2 x_0^2 \pmod{4}$, we have $y_0 \equiv p x_0 \pmod{2}$, and hence $b' \in \Z$.
By construction,
\[
x_0=u_0^{a',b'}(p,1), \quad y_0=v_0^{a',b'}(p,1).
\]
Moreover,
\[
(2b'-a'p)^2-(p^2-4)a'^2=y_0^2-(p^2-4)x_0^2 = 4N,
\]
which shows that $(a',b') \in T$.

\medskip

\noindent\textbf{Case 2:} $q = -1$.
Then $(x_0,y_0)$ satisfies
\begin{equation}\label{for q=-1}
y_0^2 - (p^2+4)x_0^2 = (-1)^n  4N
\end{equation}
for some integer $n$.

%\medskip
\noindent\underline{Subcase 2.1:} $n$ is even.
Then \eqref{for q=-1} reduces to
\[
y_0^2 - (p^2+4)x_0^2 = 4N,
\]
and the same argument as in Case~1 shows that there exist $(a',b') \in T$ such that
\[
x_0 = u_0^{a',b'}(p,-1), \quad y_0 = v_0^{a',b'}(p,-1).
\]

%\medskip
\noindent\underline{Subcase 2.2:} $n$ is odd.
Then \eqref{for q=-1} reduces to
\[
y_0^2 - (p^2+4)x_0^2 = -4N.
\]
We define $b' := x_0$ and determine $a'$ by $y_0 = b'p + 2a'$, from \eqref{Lucas seq},
which gives
\[
a' = \frac{y_0 - p x_0}{2}.
\]
Since $y_0^2 \equiv p^2 x_0^2 \pmod{4}$, we have $y_0 \equiv p x_0 \pmod{2}$, and hence $a' \in \mathbb{Z}$.
By construction,
\[
x_0 = u_1^{a',b'}(p,-1), \quad y_0 = v_1^{a',b'}(p,-1).
\]
Moreover,
\begin{align*}
(2b' - a'p)^2 - (p^2+4)a'^2
&= (2x_0 - a'p)^2 - (p^2+4)a'^2 \\
&= -\bigl(y_0^2 - (p^2+4)x_0^2\bigr) \\
&= 4N,
\end{align*}
which shows that $(a',b') \in T$.
\end{proof}

\subsection{Proof of Lemma~\ref{main theorem 2}} \label{single Pell's equations}

From Corollary~\ref{Whitney theorem} and Lemma~\ref{properties of Lucas seq1}(1),  we know that for any $n$, $(u_n^{a,b}, v_n^{a,b})$ is a solution to \eqref{gen Pell's equ}. Hence, for a fixed positive integer $p$, $q=\pm 1$, and $D=p^2-4q$, we have
\begin{equation}\label{one inclusion}
\bigcup_{(a,b)\in T}\{(\epsilon_1 u_n^{a,b}, \epsilon_2 v_n^{a,b}) \mid n\in \Z, \epsilon_1, \epsilon_2\in \{\pm 1\} \}\subset S.
\end{equation}

For the reverse inclusion, let $(x_0,y_0)\in S$. Then by Lemma~\ref{description of solution set}, $(x_0,y_0)$ is either $(u_0^{a',b'},v_0^{a',b'})$ or $(u_1^{a',b'},v_1^{a',b'})$ for some $(a',b')\in T$. Hence, $(x_0,y_0)$ belongs to the left-hand side. This completes the proof.

\subsection{Proof of Theorems~\ref{method} and \ref{method2}}

We now prove Theorems~\ref{method} and~\ref{method2}.
Let 
\begin{equation}\label{quadratic form}
Q(z,w):=w^2-pzw+qz^2
\end{equation}
 be a quadratic form and note that $(a,b)\in T$ if and only if $Q(a,b)=N$.

\paragraph{Proof of Theorem~\ref{method}}
%\noindent\textbf{Proof of Theorem~\ref{method}}
For $(a_0,b_0)\in T$ in the theorem, consider the sequence $\{u_n:=u_n^{a_0,b_0}\}_{n\in \Z}$. 
We use the following observation:
\begin{equation}\label{quadratic}
Q(u_{n},u_{n+1})=(-1)^n N.
\end{equation}
Since $Q(u_0,u_1)=u_1u_{-1}-u_0^2=N>0$ and $u_1=b_0>0$, it follows that $u_{-1}>0$.
Moreover, $u_{-2}<0$. Indeed, if $u_{-2}\geq 0$, then $Q(u_{-2},u_{-1})=N>0$, and since $u_{-2}<a_0$, this contradicts the minimality of $a_0$. We also note that $u_n\geq 0$ for $n\geq -1$.

Now let $(a',b')\in T$ and consider the Lucas sequence $\{u'_n\}_{n\in \Z}$ with $u'_m=a'$ and $u'_{m+1}=b'$ for some $m$. 
Since $Q(-z,-w)=Q(z,w)$, we may assume, after possibly changing sign and shifting the index, that $u'_{-2}<0$ and $u'_n\geq 0$ for all $n\geq -1$. 
These operations do not affect the set
\[
\{(\epsilon_1 u'_n, \epsilon_2 v'_n) \mid n\in \mathbb{Z},\ \epsilon_1,\epsilon_2\in\{\pm1\}\}.
\]

Since $Q(u'_{-2}, u'_{-1})=(u'_{-1})^2-u'_{-2}u'_0=N>0$, by \eqref{quadratic}, $Q(u'_{0}, u'_{1})=N>0$. Now the minimality of $a_0$ implies that $u'_0\geq a_0$. 

\medskip
\noindent\textit{Claim.}
Either $a_0\leq u'_0\leq b_0$ or $a_0\leq u'_{-1}\leq b_0$.

\medskip
\noindent\textit{Proof of Claim.}
If $u'_0\leq b_0$, there is nothing to show, hence, suppose $u'_0>b_0$ and we show $a_0\leq u'_{-1}\leq b_0$. 
By Lemma \ref{power of M and E}, 
\begin{equation*}\label{E-ME}
E_{(a,b)}-M^{-1}E_{(u'_0,u'_1)}=\begin{pmatrix}u_{-1}-u'_{-2}& a_0-u'_{-1}\\a_0-u'_{-1}&b_0-u'_0\end{pmatrix},
\end{equation*}
and which has negative determinant, where $M=M_{(p,-1)}$. So 
\begin{equation}
M(E_{(a,b)}-M^{-1}E_{(u'_0,u'_1)})=\begin{pmatrix} a_0-u'_{-1}&b_0-u_0'\\b_0-u_0'&b_0-u'_1\end{pmatrix}
\end{equation}
has positive determinant. Since $b_0-u_1'<0$, $a_0-u'_{-1}$ is negative, i.e., $a_0<u'_{-1}$.
Moreover, suppose $b_0< u'_{-1}$. Then since 
\begin{align*}
Q(a_0,b_0)&=b_0u_{-1}-a_0^2=N ~\ \text{ and } \\
Q(u'_{-1},u'_0)&=u'_0u'_{-2}-(u'_{-1})^2=-N,
\end{align*}
we have 
\begin{equation}\label{for contradiction}
b_0u_{-1}-a_0^2=(u'_{-1})^2-u'_0u'_{-2}.
\end{equation}
However, since $b_0u_{-1}<b_0^2<(u'_{-1})^2$ and $u'_{-2}<0$, we have a contradiction to \eqref{for contradiction}. Hence, $a_0\leq u'_{-1}\leq b_0$. \qed

By the claim, for any $(a',b')\in T$, we have an equivalence between $(a',b')$ and either $(u'_0,u'_{1})$ or $(u'_{-2},u'_{-1})$, which is in $T'$ as in \eqref{condition T'}.
Therefore, every equivalence class contains a representative in $T'$, and it suffices to consider pairs $(a,b)$ in $T'$. Therefore, by the this argument and Lemma~\ref{main theorem 2}, the set $S$ of all solutions is given by \eqref{simple form of solutions 1}.\qed

\paragraph{Proof of Theorem~\ref{method2}}
%\noindent\textbf{Proof of Theorem~\ref{method2}}
Suppose $q=1$. In this case, $Q(u_n^{a,b}, u_{n+1}^{a,b})=(u_n^{a,b})^2-u_{n+1}^{a,b} u_{n-1}^{a,b}$ is the same for all $n\in \Z$. Moreover, since $Q(u_n^{a,b}, u_{n+1}^{a,b})=(u_n^{a,b}+u_{n+1}^{a,b})^2-(p+2)u_n^{a,b}u_{n+1}^{a,b}$, $Q(a,b)<0$ if and only if $u_n^{a,b}>0$ or $u_n^{a,b}<0$ for all $n\in \Z$. 
Hence, for any $(a',b')\in T$, there is $k\in \Z$ such $u_n^{a',b'}>0$ for all $n\geq k$ and $u_n^{a,b}<0$ for all $n<k$. By the minimality of $a_0$, we have $a_0\leq u_{k}^{a',b'}$. 

Moreover, from
\begin{equation}\label{Q comparison}
\begin{aligned}
Q(a_0,b_0)&=b_0^2-pa_0b_0+a_0^2=N \\
Q(u_{k-1}^{a',b'},u_k^{a',b'})&=(u_k^{a',b'})^2-pu_{k-1}^{a',b'}u_{k}^{a',b'}+(u_{k-1}^{a',b'})^2=N,
\end{aligned}
\end{equation}
we have 
\begin{equation}
b_0^2-(u_{k}^{a',b'})^2=(pa_0b_0-a_0^2)+(u_{k-1}^{a',b'})^2-pu_{k-1}^{a',b'}u_{k}^{a',b'}>0,
\end{equation}
hence $a_0\leq u_{k}^{a',b'}<b_0$. % we may assume that $b_0\geq a_0$ because if not, then changing signs of all, we get $a_0<0<-b_0$.
Thus, any $(a',b')\in T$ is equivalent to $(u_{k}^{a',b'}, u_{k+1}^{a',b'})\in T''$ in \eqref{alpha}.

Therefore, every equivalence class contains a representative in $T''$ in \eqref{alpha}, and it suffices to consider pairs $(a,b)$ in $T''$. Therefore, by the this argument and Lemma~\ref{main theorem 2}, the set $S$ of all solutions is given by \eqref{simple form of solutions 2}.\qed

\subsection{Examples}
\begin{example}(cf. \cite[\S 5]{A})\label{example with 836}
For each $n \in \Z$, $(x,y)=(u_n^{1,15}(1,-1), v_n^{1,15}(1,-1))$ is a solution to 
\begin{equation}\label{eq in example with 836}
y^2-5x^2=\pm 836, ~\ \text{ where } ~ D=5.
\end{equation}
Although $(x,y)=(8,34)$ is also a solution to \eqref{eq in example with 836}, $8$ is not an element of $u_n^{1,\,15}(1,-1)$ for any $n\in \Z$. 

Using the method in Theorem~\ref{method}, we determine all solutions to \eqref{eq in example with 836}. For $(a,b)=(1,15)$, we obtain
\begin{equation}
T' = \{ (1,15), (5,18), (8,21), (13,27) \}.
\end{equation}
Moreover, since
\begin{equation}
\begin{aligned}
(5,18)&=(-u_{-2}^{8,21}(1,-1),u_{-3}^{8,21}(1,-1)) ~\ \text{ and } \\
(13,27)&=(-u_{-2}^{1,15}(1,-1),u_{-3}^{1,15}(1,-1)),
\end{aligned}
\end{equation}
we have $(5,18)\sim(8,21)$ and $(13,27)\sim(1,15)$.
Therefore, there are exactly two equivalence classes.

Hence, the set $S$ of all solutions to \eqref{eq in example with 836} is
\begin{equation}
\{(\epsilon_1u_n^{1,15}(1,-1),\epsilon_2v_n^{1,15}(1,-1))\} \sqcup \{(\epsilon_1u_n^{8,21}(1,-1),\epsilon_2v_n^{8,21}(1,-1))\},
\end{equation}
where $n\in \Z$ and $\epsilon_1,\epsilon_2 \in \{\pm 1\}$.
\end{example}

\begin{example}\label{example with 1276} 
For 
\begin{equation}\label{eq in example with 1276}
y^2-5x^2=\pm 1276,
\end{equation}
we may find $D=p^2-4q=5$ with $p=1, q=-1$ and $(a,b)=(2,19)$. 
Using the method given in Theorem~\ref{method}, for $(a,b)=(2,19)$, we obtain
\begin{equation}
T' = \{ (2,19),(7,23),(9,25), (15,32) \}.
\end{equation}
Since 
\begin{equation}
\begin{aligned}
(9,25)&=(-u_{-2}^{7,23}(1,-1),u_{-3}^{7,23}(1,-1)) ~\ \text{ and }\\
(15,32)&=(-u_{-2}^{2,19}(1,-1),u_{-3}^{2,19}(1,-1)),
\end{aligned}
\end{equation}
we have $(9,25)\sim(7,23)$ and $(15,32)\sim(2,19)$. 
Thus, there are exactly two equivalence classes.

Hence, the set $S$ of all solutions to \eqref{eq in example with 1276} is
\begin{equation}
\{(\epsilon_1u_n^{2,19}(1,-1),\epsilon_2v_n^{2,19}(1,-1))\}\sqcup \{(\epsilon_1u_n^{7,23}(1,-1),\epsilon_2v_n^{7,23}(1,-1))\},
\end{equation}
where $n\in \Z$ and $\epsilon_1,\epsilon_2\in \{\pm 1\}$.
\end{example}

\begin{example}\label{example with 1004} 
For 
\begin{equation}\label{eq in example with 1004}
y^2-5x^2=1004, 
\end{equation}
we may find $D=p^2-4q=5$ with $p=3, q=1$ and $(a,b)=(2,19)$.
Using the method given in Theorem~\ref{method2}, for $(a,b)=(2,19)$, we obtain
\begin{equation}
T''= \{ (2,19),(13,41) \}.
\end{equation}
Since $(13,41)=(-u_{-1}^{2,19}(3,1),-u_{-2}^{2,19}(3,1))$, the set $S$ of all solutions to \eqref{eq in example with 1004} is
\begin{equation}
\{(\epsilon_1u_n^{2,\,19}(3,1),\epsilon_2v_n^{2,\,19}(3,1))\},
\end{equation}
where $n\in \Z$ and $\epsilon_1,\epsilon_2\in \{\pm 1\}$.
\end{example}

\section{Automorphisms of K3 Surfaces}\label{Automorphisms of K3 Surfaces}

In this section, we use the properties of Lucas sequences to construct K3 surfaces equipped with automorphisms induced by the action of $M_{(p,q)}^n$ on their Picard lattices.

\subsection{K3 Surfaces}
A \textit{K3 surface} is a compact surface $X$ over $\C$ with trivial canonical bundle and $\dim H^1(X, \mathcal{O}_X)=0$ (see \cite{BHPV, H} for details). We consider the second integral cohomology $H^2(X, \Z)$ with the cup product as a lattice, that is, a free $\Z$-module of finite rank equipped with a symmetric bilinear form. It is known that $H^2(X, \Z)$ is an even unimodular lattice of signature $(3, 19)$, which is unique up to isomorphism. In the following, we denote by $S_X$, $T_X$, and $\omega_X$ the Picard lattice, the transcendental lattice and a nowhere vanishing holomorphic $2$-form on $X$, respectively. 

Let $X$ be a projective K3 surface. Since $H^2(X, \Z)$ is unimodular, we have the following natural identification:
\begin{equation}
A(S_X)\cong A(T_X)\cong H^2(X, \Z)/(S_X\oplus T_X),
\end{equation}
where, for a lattice $\Lambda$, the \textit{discriminant group} $A(\Lambda)$ of $\Lambda$ is defined as $\Lambda^*/\Lambda$, with $\Lambda^*=\rH_{\Z}(\Lambda, \Z)$ (see \cite{N} for details).

\begin{proposition}[{\cite[Proposition 2]{HKL}}] \label{Lattice of K3 rank 2}
For any even lattice $\Lambda$ of signature $(1,1)$, there exists a projective K3 surface $X$ such that $S_X\cong \Lambda$.
\end{proposition}

An embedding $S\hookrightarrow \Lambda$ of lattices is called \textit{primitive} if $\Lambda/S$ is free.
\begin{proposition}[{\cite[Corollary 1.5.2]{N}}]\label{Nikulin}
Let $S$ be a primitive sublattice of a unimodular lattice $\Lambda$, and let $\phi$ be an isometry of $S$. $\phi$ extends to an isometry of $\Lambda$ if and only if there exists an isometry $\psi$ of the orthogonal complement $T$ of $S$, such that $\bar{\psi}\circ \gamma=\gamma\circ \bar{\phi} : A(S)\xrightarrow{\sim}  A(T)$.
\end{proposition}

By the global Torelli theorem for K3 surfaces \cite{BR, PS}, an automorphism of a projective K3 surface $X$ is given by a pair of isometries $(\phi, \psi)\in O(S_X)\times O(T_X)$ such that $\phi^*|_{A(S_X)}=\psi^*|_{A(T_X)}$ and $\phi$ preserves the ample cone of $X$ and $\psi$ preserves $\C\omega_X$.

An automorphism $g$ of a K3 surface $X$ is called \textit{symplectic} if $g^*\omega_X=\omega_X$ and is called \textit{anti-symplectic} if $g^*\omega_X=-\omega_X$.

\subsection{Automorphisms of K3 Surfaces}\label{Aut of K3}
Let $X$ be a K3 surface whose Picard lattice $L\subset H^2(X,\Z)$ has signature $(1,1)$ with the intersection matrix 
\begin{equation}\label{intersection matrix of L}
\begin{pmatrix} 2r&~ s\\s &~ 2t \end{pmatrix}.
\end{equation}
We will use the notions of a lattice and its associated intersection matrix interchangeably throughout this paper.

Suppose that $M:=M_{(p,q)}$ is an induced isometry of $L$ from an automorphism $g$ of $X$.
First, note that any isometry on the Picard lattice induced by an automorphism of a K3 surface has determinant $1$, since it preserves the orientation of the lattice (cf. \cite[Proposition 3]{HKL}). Therefore, we assume that  $q\mid p$ and define $Z:=M^2/q$. 

\begin{lemma}
$Z$ is an isometry of $L$ if and only if $r=tq$, $s+tp=0$.
\end{lemma}
\begin{proof}
This follows from the fact that $Z$ is an isometry of $L$ if and only if $Z^{\transpose}LZ=L$, where $Z^{\transpose}$ is the transpose of $Z$.
\end{proof}
Now let 
\begin{equation}\label{L(t)}
L(t):=t\begin{pmatrix} 2q&~ -p\\-p &~ 2 \end{pmatrix}.
\end{equation} 

\begin{lemma}[{\cite[Lemma 1]{HKL}}]\label{criterion}
Let $\Lambda$ be an arbitrary non-degenerate even lattice of rank $n$. For $g\in O(\Lambda)$ and $\epsilon\in\{\pm 1\}$, $g$ acts on the discriminant group $A(\Lambda):=\Lambda^*/\Lambda$ of $\Lambda$ as $\epsilon\cdot \id$ if and only if $(g-\epsilon\cdot I_n)\cdot \Lambda^{-1}$ is an integer matrix, where $I_n$ is the identity matrix of size $n$. 
\end{lemma}

\begin{lemma}\label{integral conditions}
The matrix $(Z^n-I_2)L(t)^{-1}$ is integral if $tq^n\mid u_n^2$ and $tq^{n-1}\mid u_nu_{n-1}$.
\end{lemma}
\begin{proof} 
We begin by computing:
\begin{align*}
(Z^n - I_2)L(t)^{-1}
&= \frac{-1}{tq^nD}
\begin{pmatrix}
-qu_{2n-1} - q^n &~ u_{2n} \\
-qu_{2n} &~ u_{2n+1} - q^n
\end{pmatrix}
\begin{pmatrix}
2 &~ p \\
p &~ 2q
\end{pmatrix} \\
%&= \frac{-1}{tq^nD}
%\begin{pmatrix}
%pu_{2n} - 2qu_{2n-1} - 2q^n &~ q(2u_{2n} - pu_{2n-1} - pq^{n-1}) \\
%pu_{2n+1} - 2qu_{2n} - pq^n &~ q(2u_{2n+1} - pu_{2n} - 2q^n)
%\end{pmatrix} \\
&= \frac{-1}{tq^nD}
\begin{pmatrix}
u_{2n+1} - qu_{2n-1} - 2q^n &~ q(u_{2n} - qu_{2n-2} - pq^{n-1}) \\
u_{2n+2} - qu_{2n} - pq^n &~ q(u_{2n+1} - qu_{2n-1} - 2q^n)
\end{pmatrix}.
\end{align*}
Using the identities in Lemma~\ref{properties of Lucas seq1}, we compute the entries as follows:\\
The (1,1)-entry is $u_{2n+1}-qu_{2n-1}-2q^n=Du_n^2$.
The (1,2)-entry simplifies as follows:
\begin{align*}
q(u_{2n} - q u_{2n-2} - p q^{n-1})&= q\left(\frac{u_{2n+1} + q u_{2n-1}}{p}\right) - q^2 \left(\frac{u_{2n-1} + q u_{2n-3}}{p} \right) - p q^n \\
&= \frac{q}{p}(D u_n^2 + 2 q^n) + \frac{q^2}{p}(D u_{n-1}^2 + 2 q^{n-1}) - p q^n \\
&= \frac{Dq}{p} (u_n^2 + q u_{n-1}^2 - q^{n-1}) \\
&= \frac{Dq}{p} u_{n-1} (q u_{n-1} + u_{n+1}) \\
&= D q u_n u_{n-1}.
\end{align*}
Similarly, the (2,1)-entry becomes $Dqu_nu_{n+1}$ and the $(2,2)$-entry is $Dqu_n^2$. 

Thus,
\begin{equation}\label{power-identity}
(Z^n-I_2)L(t)^{-1}=\frac{-u_n}{tq^n}\begin{pmatrix} u_n&~ qu_{n-1}\\ u_{n+1} &~ qu_n \end{pmatrix}.
\end{equation}
Since $u_{n+1}=pu_n-qu_{n-1}$, the above matrix has integer entries provided that $tq^n\mid u_n^2$ and $tq^{n-1}\mid u_nu_{n-1}$.
\end{proof}

Let $g$ be an automorphism of a K3 surface $X$. For the induced isometry $g^*:H^2(X,\C)\rightarrow H^2(X,\C)$, let $\lambda(g)$ be the spectral radius:
\begin{equation}\label{spectral radius}
\lambda(g):=\max\{\vert \mu\vert \mid \mu\in \C \text{ is an eigenvalue of } g^*\}.
\end{equation}
It is known that $\lambda(g)$ is equal to $1$ or a \textit{Salem number}, i.e., a real algebraic number $\lambda>1$ whose conjugates other than $\lambda^{\pm 1}$ lie on the unit circle (see \cite{M,M2} for details). 
%By the fundamental result of Gromov and Yomdin ([\cite{G,Y}]), the topological entropy of an automorphism of a compact K\"ahler surface is the logarithm of a Salem number. 
The irreducible monic integer polynomial $S(x)$ of $\lambda$ is a \emph{Salem polynomial} and the degree of $S(x)$ is the degree of Salem number $\lambda$. A \emph{Salem trace} is an algebraic integer $\tau>2$ whose other conjugates lie in $(-2,2)$. Salem traces and Salem numbers correspond bijectively, via the relation $\tau=\lambda+\lambda^{-1}$.

\begin{theorem}\label{K3 automorphisms from matrix representations}
Suppose $q \neq \pm 1$, $q \mid p$, $D=p^2-4q>0$, and assume that $D$ is not a square number and $L:= L(1)$ contains no $(-2)$ vectors, i.e., the lattice $L$ has no self-intersection $-2$ element. Then there exists a K3 surface $X$ with Picard lattice $L$ admitting a symplectic automorphism $g$ such that
\begin{equation}
g^*|_L=Z^2.
\end{equation}
In particular, for any $n \in \Z$, the Salem trace of $g^n$ is given by
\begin{equation}
\tr \left(Z^{2n}\right) = \left( \frac{v_{2n}}{q^n} \right)^2-2.
\end{equation}
\end{theorem}

\begin{proof}
By Proposition~\ref{Lattice of K3 rank 2}, there exists a projective K3 surface $X$ whose Picard lattice is $L$.
Since $D$ is not a square, the lattice $L$ contains no elements of square zero. Indeed, for an element $(x,y)\in L$, the self-intersection number is given by $2(y^2-pxy+qx^2)$, which is nonzero for all $(x,y)\neq (0,0)$. 

Since $q \mid p$ and $u_2=p$, Lemmas~\ref{criterion} and~\ref{integral conditions} imply that the matrix $Z^{2}$ is an isometry of $L$ and acts trivially on the discriminant group $A(L) := L^*/L$. Therefore, it extends to an isometry of $H^2(X, \mathbb{Z})$ by acting as the identity on the transcendental lattice $T_X$ by Proposition \ref{Nikulin}.

Since $L$ contains no $(-2)$-vectors, the ample cone coincides with the positive cone (cf. \cite[Corollary 8.2.11]{H}). Because $Z^{2}$ preserves the ample cone, the global Torelli theorem implies that it lifts to an automorphism $g$ of the K3 surface $X_L$ with $g^*|_L = Z^{2}$.

Finally, we compute the Salem trace; from Corollary~\ref{Whitney theorem} and Lemma~\ref{properties of Lucas seq1}(1), we obtain
\begin{equation*}
\tr(g^{*n}|_L)=\tr Z^{2n} =\frac{\tr(M^{4n})}{q^{2n}}=\frac{Du_{2n}^2+4q^{2n}}{q^{2n}}-2=\left( \frac{v_{2n}}{q^n} \right)^2-2.
\end{equation*}
\end{proof}

\begin{remark} 
For $q=-1$, we refer Theorem 1.1 in \cite{L1}.
\end{remark}

\section*{Acknowledgments}
This research was supported by the Basic Science Research Program through the National Research Foundation of Korea(NRF), funded by the Ministry of Education (2019R1A6A1A10073887); by the National Research Foundation of Korea(NRF) grant funded by the Korea government(MSIT) (RS-2024-00352592); and by the Institute for Basic Science (IBS-R032-D1).

\end{document}